\RequirePackage{fix-cm}
	\documentclass[twocolumn]{svjour3}          
\smartqed  
\usepackage{graphicx}
\usepackage{amssymb}
\usepackage{amsmath,color,latexsym} 
\newcommand{\eps}{\varepsilon}
\usepackage{extarrows}

 \journalname{}

\newcommand{\sign}{{\rm \hskip0.5pt sign \hskip1pt}}

\begin{document}

\title{A new test for stick-slip limit cycles in dry-friction oscillators with small nonlinear friction characteristics 
}


\author{        Oleg Makarenkov 
}


\institute{           O. Makarenkov \at
              Department of Mathematical Sciences, University of Texas at Dallas, 75080 Richardson, USA\\
              \email{makarenkov@utdallas.edu}
}

\date{Received: date / Accepted: date}

\maketitle

\begin{abstract} We consider a dry friction oscillator on a moving belt
with both the Coulomb friction and a small nonlinear addition which can model e.g. the Stribeck effect. By using the perturbation theory, we establish a new sufficient condition for the nonlinearity to ensure the occurrence of a stick-slip limit cycle. The test obtained is more accurate compared to what one gets by building upon the divergence test.

\end{abstract}

\section{Introduction} 
\begin{figure}[h]\center
\vskip-0.3cm
\includegraphics[scale=0.77]{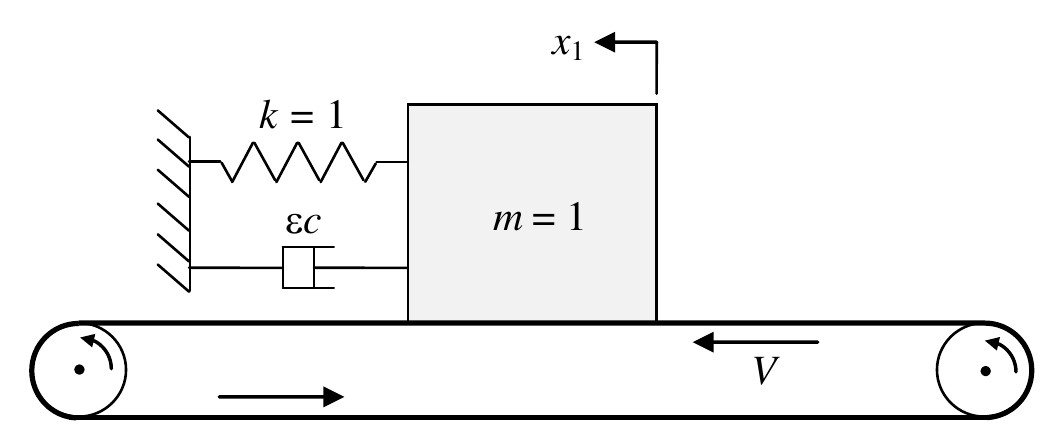}
\caption{\footnotesize Dry friction oscillator on a moving belt.} \label{mass_belt}
\end{figure}

Dry friction oscillator on a moving belt (Fig.~\ref{mass_belt}) is a fundamental model of classical mechanics. One of the main aims of this model is its ability to resemble the stick-slip oscillations crucial in understanding the tectonic fault in seismology \cite{earth1,earth2,earth3}, squealing noise in automotive industry \cite{sq1,sq2,sq3}, hysteresis effect in drilling systems \cite{drill}. Being a simplified model, the dry friction oscillator is not capable to accurately simulate the dynamics of a tectonic fault or a brake squeal. The role of simplified (conceptual) models is different. They are rather supposed to provide handy maps that unveil important phenomena and computable conditions for the parameters to be further tested in experiments or in simulations on comprehensive models (see an interesting discussion in \cite{climate}).  Indeed, derivation of suitable Poincar\'e maps allowed the authors of \cite{szalai,kry} to prove the existence of chaotic sticking motions in the model of Fig.~\ref{mass_belt}. The work \cite{gal} uses the divergence analysis to link the lack of stick-slip oscillations to the parameters of Stribeck friction characteristics. The authors of \cite{closed1,closed2,pascal} derive conditions for stick-slip oscillations in particular equations of mass-spring oscillator on a frictional belt by computing  closed-form equations of motion.

\noindent 
In this paper we offer sufficient conditions for the occurrence of stick-slip oscillations in the model of Fig.~\ref{mass_belt} assuming the friction characteristics is just arbitrary, but small. Specifically, we assume that the oscillator under consideration is modeled by 
\begin{equation}\label{dfo}
\begin{array}{l}
   \dot x_1=x_2,\\
   \dot x_2+\eps c x_1 +x_2=-{\rm sign}(x_2\hskip-0.05cm-\hskip-0.05cm V)\left(1\hskip-0.05cm+\hskip-0.05cm F(x_2\hskip-0.05cm-\hskip-0.05cm V,\eps)\right),
\end{array}
\end{equation}
where $\eps>0$ is a small parameter and 
where the dry friction characteristics is a combination of the Coulomb term ${\rm sign}(x_2-V)$ and an arbitrary  Lipschitz ingredient $F(x_2-V,\eps)$ which satisfies
\begin{equation}\label{F}
   F(0,\eps)\equiv F(x_2,0)\equiv0.
\end{equation}
 The main result of the paper states that stick-slip oscillations take place, if 
\begin{equation}\label{achiev}
  cV\pi<\int_0^{2\pi}F'_\eps(V\cos\tau-V,0)\cos\tau d\tau.
\end{equation}
The method of proof is perturbation theory that already demonstrated its efficiency in various dry friction problems, see \cite{fidlin,awr}.

\vskip0.2cm

\noindent The paper is organized as follows. The next section spots the line segment of the line $x_2=V$ where the solutions of (\ref{dfo}) stick to $x_2=V$. To discuss sticking rigorously, the Filippov definition of solution of discontinuous differential equations is introduced and some basic background is presented. In section 3 we show that the role of nonlinearity in the occurrence of stick-slip oscillations in (\ref{dfo}) is pivotal. Specifically, we show that just Coulomb friction along is insufficient to produce stick-slip limit cycles, if the viscous friction is small. The main result of the paper (Proposition~\ref{pr2}) is proved in section 4,  where we derive a condition to ensure that the initial condition $(1,V)$ whose $2\pi$-periodic trajectory just touches the line $x_2=V$ when $\eps=0$ transforms into an initial condition $(1-\eps c V,V)$ those trajectory does return back to $x_2=V$ over the time $2\pi+O(\sqrt{\eps})$ and slides along $x_2=V$ until it reaches $(1-\eps c V,V)$. Such a transition is know as grazing-sliding bifurcation in the field of nonsmooth dynamical systems, see \cite{lamb}. In section 5 we specify the main result of the paper for the case of Stribeck nonlinearity and identify the parameters of Stribeck nonlinearity (Proposition~\ref{prstibeck}) which lead to the occurrence of a stick-slip limit cycle. In section 5.1 we compare Proposition~\ref{prstibeck} with what one obtains by using the divergence based test.

\vskip0.2cm 

\noindent We would like to stress that the focus of this work is on the development of a rigorous mathematical theory of the dry-friction oscillator. Readers interested in simulation results are referred to other papers extensively available in the literature, see
e.g. \cite{wie,sq1,gal,kuepper}.

\section{Sliding along the switching threshold}

We remind the reader that for a system of piecewise smooth differential equations
$$
  \dot x =f(x)
$$
discontinuous along the line $x_2=V,$ a Filippov solution $x$ is an absolutely continuous function that satisfies
$$
   \dot x(t)\in K[f](x(t)),\quad{\rm for\ a.\ a\ }t,
$$
where $K[f](x)=\bigcup\limits_{\lambda\in[0,1]}(\lambda f(x,V^-)+(1-\lambda)f(x,V^+))$ is the convexification of the discontinuous function $f$ and
\begin{eqnarray*}
f(x_1,V^-)&=&\lim\limits_{x_2\to V-0} f_2(x_1,x_2),\\ 
f(x_1,V^+)&=&\lim\limits_{x_2\to V+0} f_2(x_1,x_2).
\end{eqnarray*}
The
 point $(x_1,V)$ is a point of sliding, if
$$
\lim\limits_{x_2\to V^-} f_2(x_1,x_2)>0\quad{\rm and}\quad \lim\limits_{x_2\to V^+} f_2(x_1,x_2)<0.
$$
Furthermore, the Filippov equation of sliding along $x_2=V$ is 
$$
   \dot x\in K[f](x)\cap \mathbb{R}\times\{0\},
$$
because the set of all vectors tangent to the switching manifold $x_2=V$ is $\mathbb{R}\times\{0\}$  at any point of this switching manifold.

\vskip0.2cm

\begin{figure}[h]\center
\includegraphics[scale=0.77]{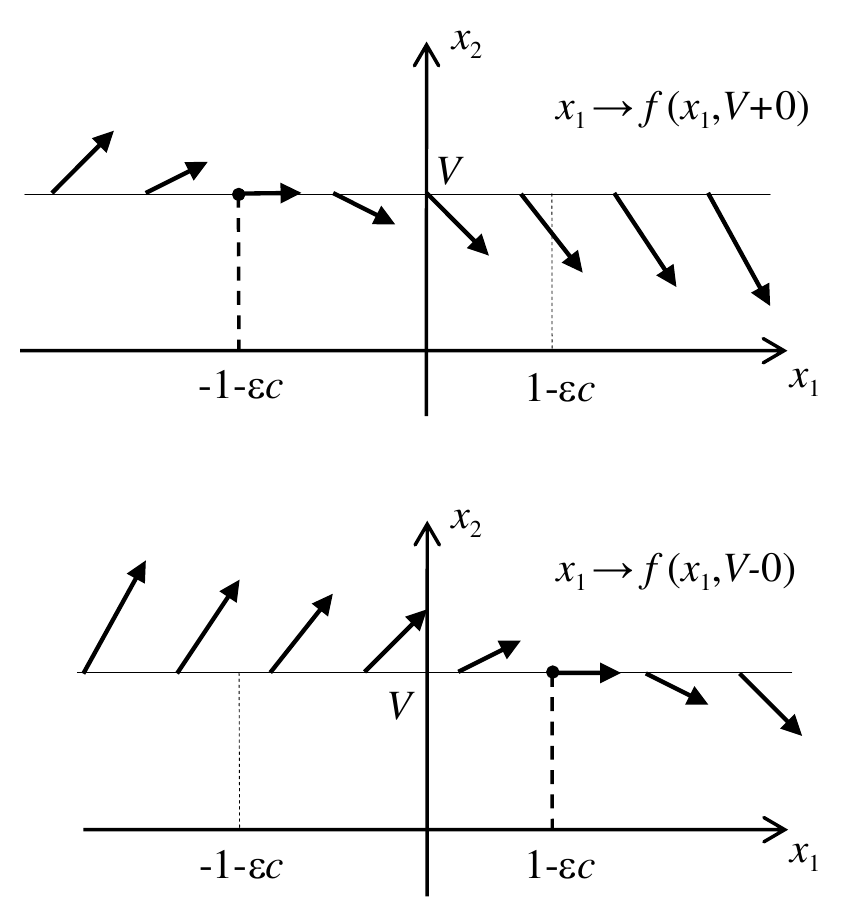}
\caption{\footnotesize The behavior of vector fields $x_1\to f(x_1,V^-)$ and $x_1\to f(x_1,V^+)$  of (\ref{dfo}) along the line $x_2=V$.} \label{dry_friction}
\end{figure}

\noindent Consequently, a point $(x_1,V)$ is a point of sliding for (\ref{dfo}), if
$$
   -x_1-\eps cV+1>0\quad{\rm and}\quad -x-\eps cV-1<0,
$$
see Fig.~\ref{dry_friction}.
The Filippov equation of sliding along $x_2=V$ becomes
$$
\begin{array}{l}
   \dot x_1=V,\\
   \dot x_2=0,
\end{array}
$$
as far as system (\ref{dfo}) is concerned. The following simple proposition will be our main tool in spotting stick-slip limit cycles for (\ref{dfo}).

\begin{proposition}\label{ppp}  Let $(x_1(t),x_2(t))$ be a Filippov solution of (\ref{dfo}) with the initial condition $(x_1(0),x_2(0))=(1-\eps c,V)$. If 
$$
   x(\bar t)\in  (-1-\eps c V,1-\eps c V)\quad \mbox{and}\quad x_2(\bar t)=V
$$ 
for some $\bar t> 0,$ then $(x_1(t),x_2(t))$ is a finite-time stable limit cycle of (\ref{dfo}).
\end{proposition}

\noindent The proof is a combination of 3 facts: 1) any solution of (\ref{dfo}) with the initial condition in $(-1-\eps cV,1-\eps cV)\times \{V\}$ reaches the point $(1-\eps cV,V)$;
2) without loss of generality we can assume that $x_2(t)<V,$ $t\in(0,\bar t)$, and $(x_1(t),x_2(t))$ is the unique solution of (\ref{dfo}) on the interval $(0,\bar t)$;
3) the solution $(x_1(t),x_2(t))$ reaches $x_2=V$ at $t=\bar t$ transversally and so the solutions of (\ref{dfo}) that originate near $(x_1(0),x_2(0))$ do reach $x_2=V$ too.

\section{The case of Coulomb dry friction} In this case $f\equiv 0$ and equation (\ref{dfo}) takes the form
\begin{eqnarray}
& & \begin{array}{l}
   \dot x_1=x_2,\\
   \dot x_2=-x_1-\eps cx_2-1,\label{col1}
\end{array}\qquad\mbox{if} \ x_2>V,\\
& & \begin{array}{l}
   \dot x_1=x_2,\\
   \dot x_2=-x_1-\eps cx_2+1,
\end{array}\qquad\mbox{if} \ x_2<V.\label{col2}
\end{eqnarray}

\begin{figure}[h]\center
\vskip-0.5cm
\includegraphics[scale=0.77]{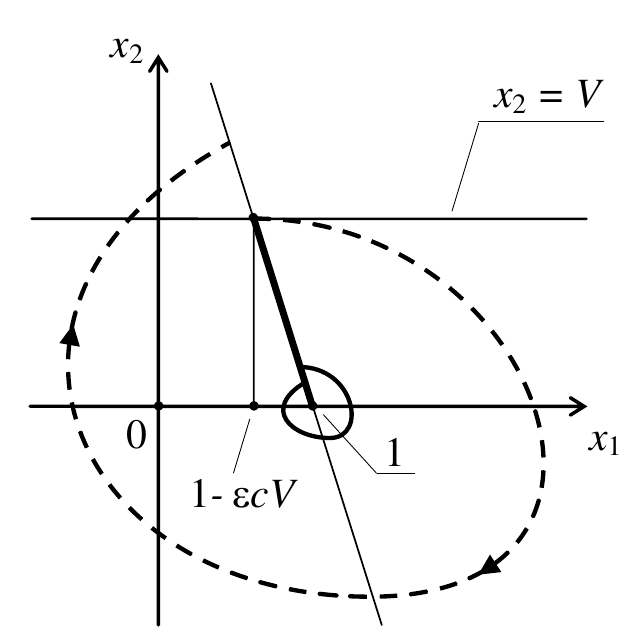}
\caption{\footnotesize True (solid curve) and false (dashed curve) trajectories of (\ref{col1})-(\ref{col2}).} \label{dfp}
\end{figure}

\noindent The system (\ref{col2}) possesses a globally stable equilibrium $(1,0)$. Let us  show that this implies that the solution $(x_1(t),x_2(t))$ of (\ref{col2}) with the initial condition $(x_1(0),x_2(0))=(1-\eps c,V)$ intersects the $x_2>0$ part of the cross-section $-x_1-\eps cx_2+1=0$ at those $t>0$ for which $x_2(t)\in (0,V)$ (i.e. the solutions that originate in the bold line segment of Fig.~\ref{dfp} land into that segment again). Indeed, assume the contrary, i.e. assume that $(x_1(t),x_2(t))$ intersects the $x_2>0$ part of the line $-x_1-\eps cx_2+1=0$ at those $t>0$ for which $x_2(t)>V$ (dashed curve in Fig.~\ref{dfp}). But, since $(1,0)$ is asymptotically stable, we can construct a second trajectory that spirals towards $(1,0)$ (solid curve in Fig.~\ref{dfp}) and conclude the existence of a limit cycle for (\ref{col2}), which cannot happen because $(1,0)$ is a globally asymptotically stable equilibrium of (\ref{col2}). The contradiction obtained proves that system (\ref{col1})-(\ref{col2}) doesn't have closed orbits passing through the sliding region $(-1-\eps cV,1-\eps cV)\times \{V\}$.

\section{The main result}

In this section we prove a general result as for how introducing a nonlinearity makes the existence of stick-slip limit cycles possible.

\begin{proposition}\label{pr2} Let $(x_2,\eps)\mapsto F(x_2,\eps)$ be a Lipschitz function which is differentiable in $\eps$ and satisfies (\ref{F}). If
\begin{equation}\label{condi}
cV\pi<\int\limits_0^{2\pi} F'_\eps(V\cos\tau -V,0)\cos\tau \hskip0.05cm d\tau,
\end{equation}
then, for all $\eps>0$ sufficiently small, the Filippov solution $x$ of (\ref{dfo}) with the initial condition $x(0)=(1-\eps c V,V)^T$ reaches the point $(\xi_1,V)^T$, where $\xi_1\in (-1-\eps c V,1-\eps c V).$ As a consequence, solution $x$ is a finite-time stable limit cycle of (\ref{dfo}).
\end{proposition}

\noindent {\bf Proof.} When $x_2<V$, system (\ref{dfo}) takes the form
\begin{equation}\label{replace1}
 \dot x=g(x,\eps),
\end{equation}
where
$$
 g(x,\eps)=
   \left( \begin{array}{llll}
       x_2 \\
       -x_1-\eps cx_2+1+F(x_2-V,\eps) 
    \end{array}\right).$$
Let $t\mapsto X(t,\xi,\eps)$ be the solution $x$ of (\ref{replace1}) with the initial condition $x(0)=\xi.$ To prove the proposition, we have to establish the existence of a function $T(\eps)$ such that 
$$
X_1\left(T(\eps),\left(\begin{array}{c} 1-\eps c V \\ V\end{array}\right),\eps\right)<1-\eps c V,$$
$$
X_2\left(T(\eps),\left(\begin{array}{c} 1-\eps c V \\ V\end{array}\right),\eps\right)=V
$$
and for which 
$$
X_2\left(t,\left(\begin{array}{c} 1-\eps c V \\ V\end{array}\right),\eps\right)<V\quad\mbox{for all}\ \ t\in(0,T(\eps)).
$$
Since the trajectory $t\mapsto X\left(t,\left(\begin{array}{c} 1-\eps c V \\ V\end{array}\right),\eps\right)$ approaches the circle of period $2\pi$ as $\eps\to 0$, we look for the solution $T(\eps)$ that approaches $2\pi$ as $\eps\to 0.$

\vskip0.3cm

\begin{figure}[t]\center
\includegraphics[scale=0.77]{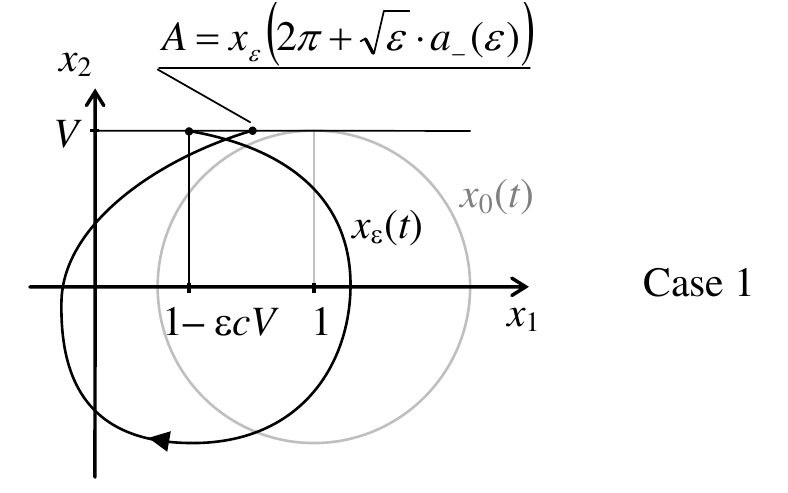}
\includegraphics[scale=0.77]{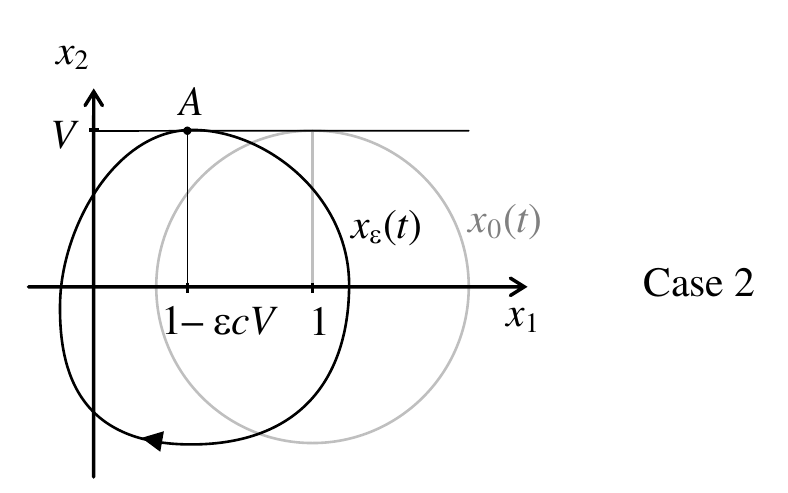} \vskip0.2cm
\includegraphics[scale=0.77]{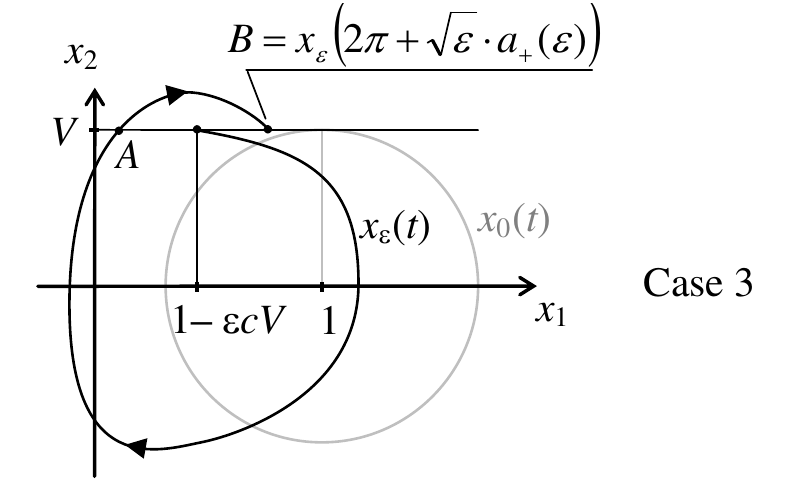}
\caption{\footnotesize Three possible locations of solution $x_\eps$ with respect to $x_0$.} \label{cases}
\end{figure}

\noindent {\bf Step 1.} Consider\footnote{Similar analysis can be execute by letting $\eps=a\delta^2$ and by considering $$
G(\delta,a)=\frac{1}{\delta}\left(X_2\left(2\pi-\delta,\left(\begin{array}{l} 1-a\delta^2 cV\\ V\end{array}\right),a\delta^2\right)-V\right).
$$} 
$$
  G(a,\eps)=\frac{1}{\eps}\left(X_2\left(2\pi+\sqrt{\eps}\cdot a,\left(\begin{array}{c} 1-\eps c V \\ V\end{array}\right),\eps\right)-V\right),
$$
which can be viewed as a {\it desingularization} or {\it blow-up} of the function $\eps G(a,\eps)$ at $\eps=0$.
To apply the Implicit Function Theorem (see e.g. \cite[\S~8.5.4, Theorem~1]{zorich}), we need to find $a\in\mathbb{R}$ such that 
$$
    G(a,0)=0,
$$
where $G(a,0)=\lim_{\eps\to 0}G(a,\eps)$ by definition. In what follows,  the abbreviation $"\xlongequal{L}"$ stays for "by applying the l'Hospital rule". Letting
$$
  \phi(a,\eps)=2\pi+\sqrt{\eps}\cdot a,\left(\begin{array}{c} 1-\eps c V \\ V\end{array}\right)
$$
we have
\begin{eqnarray*}
  G(a,0) &\xlongequal{L}&\lim\limits_{\eps\to 0}\hskip0.1cm\frac{\partial}{\partial\eps}\left(X_2(\phi(a,\eps),\eps)-V\right)=
\end{eqnarray*}
\begin{eqnarray*}
&=&\lim\limits_{\eps\to 0}\left(\frac{X_2{}'_t\left(\phi(a,\eps),\eps\right)}{2\sqrt{\eps}}a+\right.\\
&&+\left.X_2{}'_{x_1}\left(\phi(a,\eps),\eps\right)\cdot (-cV)+X_2{}'_\eps\left(\phi(a,\eps),\eps\right)\right)=
\end{eqnarray*}
\begin{eqnarray*}
&\xlongequal{L}&\frac{1}{2}X_2{}'_t{}'_t\left(2\pi,\left(\begin{array}{c} 1 \\ V\end{array}\right),0\right)a^2-\\
&&-cV \cdot X_2{}'_{x_1}\left(2\pi,\left(\begin{array}{c} 1 \\ V\end{array}\right),0\right)+X_2{}'_\eps\left(2\pi,\left(\begin{array}{c} 1 \\ V\end{array}\right),0\right).
\end{eqnarray*}
Since the property
$$
  X_2\left(2\pi,\left(\begin{array}{c} \xi_1 \\ V\end{array}\right),0\right)=V\quad\mbox{for any}\ \xi_1\in\mathbb{R}
$$
implies 
$$ X_2{}'_{x_1}\left(2\pi,\left(\begin{array}{c} 1 \\ V\end{array}\right),0\right)=0,
$$
the equation
$$
   G(a,0)=0
$$
yields
\begin{equation}\label{sqrt}
   a=\pm \sqrt{-2\cdot \frac{X_2{}'_\eps\left(2\pi,\left(\begin{array}{c} 1 \\ V\end{array}\right),0\right)}{X_2{}'_t{}'_t\left(2\pi,\left(\begin{array}{c} 1 \\ V\end{array}\right),0\right)}}.
\end{equation}
The two roots correspond to the two intersections of the solution
$$
   x_\eps(t)=X\left(t,\left(\begin{array}{c} 1 -\eps cV\\ V\end{array}\right),\eps\right)
$$
with the line $x_2=V,$ see Fig.~\ref{cases} (case 3).

\vskip0.3cm

\noindent To check the condition
$$
   G'_a(a,0)\not=0
$$
of the Implicit Function Theorem, we compute
$$
   G'_a(a,\eps)=\frac{X_2{}'_t\left(2\pi+\sqrt{\eps}\cdot a,\left(\begin{array}{c} 1-\eps c V \\ V\end{array}\right),\eps\right)}{\sqrt{\eps}}.
$$
Therefore,
$$
   G'_a(a,0) \xlongequal{L}X_2{}'_t{}'_t \left(2\pi,\left(\begin{array}{c} 1 \\ V\end{array}\right),0\right).
$$

\noindent {\bf Step 2.} Computing $X_2{}'_t{}'_t \left(2\pi,\left(\begin{array}{c} 1 \\ V\end{array}\right),0\right)$. By taking the derivative of 
\begin{equation}\label{td}
X{}'_t \left(t,\left( 1, V\right)^T,\eps\right)=g\left(X \left(t,\left( 1, V\right)^T,\eps\right),\eps\right),
\end{equation}
in $t$, 
one gets
\begin{eqnarray*}
&& X_2{}'_t{}'_t \left(2\pi,\left( 1, V\right)^T,0\right)=\\
&&\hskip-0.2cm=g_2{}'_x\left(X \left(2\pi,\left( 1, V\right)^T,0\right),0\right)X'_t\left(2\pi,\left( 1, V\right)^T,0\right)=\\
&&\hskip-0.2cm=g_2{}'_x\left(X \hskip-0.1cm\left(2\pi,\left( 1, V\right)^T\hskip-0.1cm,0\right),0\right)g\left(X\hskip-0.1cm \left(2\pi,\left( 1, V\right)^T\hskip-0.1cm,0\right),0\right).
\end{eqnarray*}
The function $t\mapsto X \left(t,\left( 1, V\right)^T,0\right)$ is the solution of the linear system
\begin{equation}\label{linsys}
  \begin{array}{lll}
     \dot x_1& = &x_2\\
     \dot x_2&=&-x_1 +1.
   \end{array}
\end{equation}
Thus, $X \left(2\pi,\left( 1, V\right)^T,0\right)=X \left(0,\left( 1, V\right)^T,0\right)=\left( 1, V\right)^T$ and
$
    X_2{}'_t{}'_t \left(2\pi,\left(\begin{array}{l}
1\\ V
 \end{array}\right),0\right)= (-1,0)\left(\begin{array}{l}
V\\ 0
 \end{array}\right)=-V.
$

\noindent{\bf Step 3.}  Computing $X_2{}'_\eps \left(2\pi,\left(\begin{array}{c} 1 \\ V\end{array}\right),0\right)$. Let $y(t)=X{}'_\eps \left(t,\left( 1, V\right)^T,0\right).$ 
By taking the derivative of (\ref{td}) in $\eps$ one gets
$$
  \dot y=g{}'_x\left(X \left(t,\left( 1, V\right)^T,0\right),0\right)y+g{}'_\eps\left(X \left(t,\left( 1, V\right)^T,0\right),0\right)
$$
or
\begin{equation}\label{lp}
\begin{array}{l}
  \dot y=\left(\begin{array}{ll} 0 & 1 \\ -1 & 0 \end{array}\right)y\hskip0.1cm+\\
\left(\begin{array}{l} 0 \\ \hskip-0.1cm-c X_2 (t,( 1, V)^T,0)\hskip-0.1cm+\hskip-0.1cm F'_\eps\left(X_2(t,( 1, V)^T,0)\hskip-0.1cm-\hskip-0.1cmV\right)\end{array}\hskip-0.1cm\right).
\end{array}
\end{equation}
By solving (\ref{linsys}), 
$$
  X_2\left(t,\left( 1, V\right)^T,0\right)=V\cos(t),
$$
and the method of variation of constants yields
$$
\begin{array}{l}
  y(t)=Z(t)y(0)+\\
+Z(t)\int\limits_0^t Z(-\tau) \left(\begin{array}{l} 0 \\ -c V\cos\tau+F'_\eps\left(V\cos \tau-V\right)\end{array}\right) d\tau,
\end{array}
$$
where $Z(t)$ is the solution of the linear part of (\ref{lp}) given by
$$
   Z(t)=\left(\begin{array}{ll} \cos t & \sin t\\  -\sin t & \cos t\end{array}\right).
$$
By noticing that $y(0)=0$, we finally conclude
\begin{eqnarray*}
   y_2(2\pi)&=&\int\limits_0^{2\pi}
\cos\tau (-cV\cos\tau+F'_\eps\left(V\cos\tau-V,0\right))d\tau=\\
&=&-cV\pi+\int\limits_0^{2\pi}F'_\eps(V\cos\tau -V,0)\cos\tau \hskip0.05cm d\tau.
\end{eqnarray*}

\noindent {\bf Step 4.} Substituting the results of Step 2 and Step 3 into (\ref{sqrt}), and by using assumption (\ref{condi}), we conclude the existence of  $a_-<0$ and $a_+>0$ such that 
$$
   G(a_-,0)=G(a_+,0)=0\ \mbox{and}\ G'_a(a_-,0)\cdot G'_a(a_+,0)\not=0.
$$
Therefore, by the Implicit Function Theorem, there exist
$$
    a_-(\eps)\to a_-\quad\mbox{and}\quad a_+(\eps)\to a_+\quad\mbox{as}\quad \eps\to 0
$$
such that 
$$
  G(a_-(\eps),\eps)=G(a_+(\eps),\eps)=0
$$
for all $\eps>0$ sufficiently small. This implies
$$
  x_{2,\eps}\left(2\pi+\sqrt{\eps}a_-(\eps)\right)=x_{2,\eps}\left(2\pi+\sqrt{\eps}a_+(\eps)\right)=V
$$
for all $\eps>0$ sufficiently small.
The following three cases are possible now (see Fig.~\ref{cases}).
\vskip0.2cm

\noindent\underline{Case 1:} $x_{1,\eps}\left(2\pi+\sqrt{\eps}a_-(\eps)\right)>1-\eps c V.$ This case is impossible because it implies that the solution $x_\eps$ crosses itself on $(0,2\pi+\sqrt{\eps}a_-(\eps))$, which cannot happen because of uniqueness of solutions of (\ref{replace1}).
\vskip0.2cm

\noindent\underline{Case 2:} $x_{1,\eps}\left(2\pi+\sqrt{\eps}a_-(\eps)\right)=1-\eps c V.$  In this case $x_\eps$ is a $\left(2\pi+\sqrt{\eps}a_-(\eps)\right)$-periodic solution of (\ref{replace1}) that intersects $x_2=V$ at only one point $(1-\eps c V,V)$. This contradicts the existence of the second intersection 

\noindent \centerline{$\left(x_{1,\eps}\left(2\pi+\sqrt{\eps}a_+(\eps)\right),V\right).$}

\vskip0.2cm

\noindent\underline{Case 3:} $x_{1,\eps}\left(2\pi+\sqrt{\eps}a_-(\eps)\right)<1-\eps c V.$ This is what was required to prove.

\vskip0.2cm

\noindent The proof of the proposition is complete.\qed

\section{The case of Stribeck nonlinearity} 
Stribeck effect is a particular type of nonlinearity in dry friction characteristics that exhibits non-monotonicity with respect to velocity (see \cite[\S 4.1-4.2]{leine}). The following function
$$
   F(x_2,\eps)=\dfrac{1-\alpha}{1+\eps\gamma|x_2|}+\alpha+\eps\beta x_2^2-1,
$$
is one of the common ways 
to introduce Stribeck effect in the dry friction oscillator of Fig.~\ref{mass_belt}, see \cite{sq1,gal,nature,kuepper} (where graphs of Stribeck dry friction for various parameters are also available). Here $\alpha,\beta,\gamma>0$ are fixed constants and $\eps>0$ measures how strong the Stribeck effect is.

\vskip0.2cm

\noindent In this section we apply Proposition~\ref{pr2} in order to determine the constants $\alpha,\beta,\gamma$ that ensure the occurrence of stick-slip oscillations when $\eps>0$ crosses zero.  To this end we compute $F'_\eps$ and get
$$
  F'_\eps(x_2,0)=\gamma(1-\lambda)x_2+\beta x_2^2.
$$
Accordingly, (\ref{condi}) takes the form
$$
   cV\pi<-\pi V(\alpha \gamma+2\beta V-\gamma)
$$
or, equivalently,
\begin{equation}\label{abc}
  -c+\gamma(1-\alpha)-2\beta V>0.
\end{equation}
To summarize, the following results is obtained. 

\begin{proposition}\label{prstibeck}
If (\ref{abc}) holds then, for all $\eps>0$ sufficiently small, the dry friction oscillator with Stribeck effect
\begin{equation}\label{dfo1}
\begin{array}{l}
\dot x_1=x_2,\\
\dot x_2=-x_1-\eps c x_2
-\sign(x_2-V)\cdot\\ 
\qquad\cdot\left(\dfrac{1-\alpha}{1+\eps\gamma|x_2-V|}+\alpha+\eps\beta(x_2-V)^2\right)
\end{array}
\end{equation}
admits a finite-time stable stick-slip limit cycle $x_\eps$ that passes through the point $(1-\eps c V,V)$.
\end{proposition}

\subsection{Comparison with the divergence test} 

In this section we offer an alternative simpler approach to derive sufficient conditions for the existence of stick-slip limit cycles in (\ref{dfo1}). Following, \cite{llibre}, this approach is based on computing the divergence.  In the context of the model (\ref{dfo1}) the divergence is used in \cite{gal} to prove the lack of limit cycles.  

\vskip0.2cm

\begin{figure}[h]\center
\includegraphics[scale=0.77]{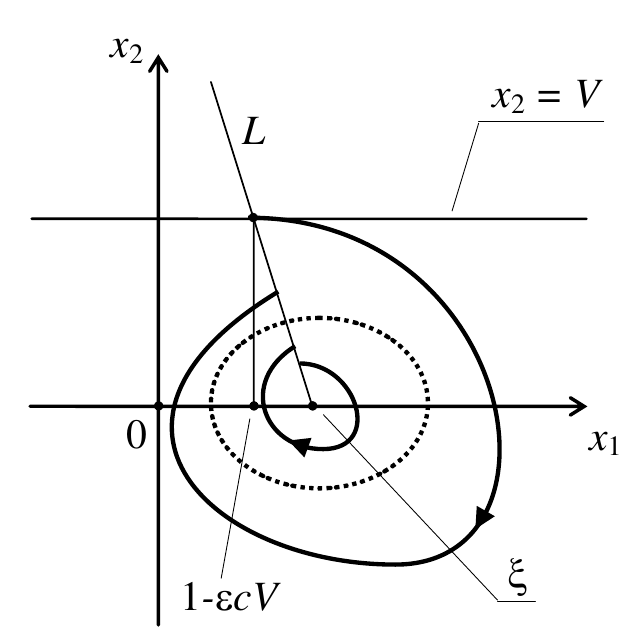}
\caption{\footnotesize Illustration of the proof of the lack of limit cycles in (\ref{ns}) which live entirely below $x_2=V.$} \label{nfig}
\end{figure}

\noindent The equilibrium of this system is $(\xi,0)=\left(\dfrac{1-\alpha}{1+\eps\gamma V}+\alpha+\eps\beta V^2,0\right)$, which is unstable if
\begin{equation}\label{unstab}
   -c+\dfrac{\gamma(1-\alpha)}{(1+\eps\gamma V)^2}-2\beta V>0.
\end{equation}
Note, that for all $\eps>0$ sufficiently small, (\ref{unstab}) is a consequence of (\ref{abc}). 
We assume that (\ref{unstab}) holds from now on.
Denote by $L$ the $x_2>0$ part of the line through $(\xi,0)$ and $(1-\eps cV,V)$, see Fig.~\ref{nfig}. If $\eps>0$ is small enough, then the solution $(x_1(t),x_2(t))$ of equation (\ref{dfo1}) with the initial condition $(x_1(0),x_2(0))=(1-\eps cV,V)$ will intersect the $x_2>0$ part of line $L$ at some time moment $\tau>0$ again. Assume that $\tau>0$ is the first time moment when $(x_1(t),x_2(t))$ intersects the $x_2>0$ part of $L$. There are two cases to consider:

\vskip0.2cm

\noindent {\boldmath $x_2(\tau)\le 1:$} This solution $x$ is the bigger solid curve in Fig~\ref{nfig}.
Since the equilibrium $(\xi,0)$ is unstable, then there is an orbit $(\hat x_1(t),\hat x_2 (t))$ with the initial condition on $L$ and close to $(x_0,0)$ such that for the next intersection $(\hat x(\hat t),\hat y(\hat t))$ of this orbit with $L$ we have $\hat y(\hat t)>\hat y(0)$  (smaller dotted curve in Fig~\ref{nfig}). 
This implies that the system
\begin{equation}\label{ns}
\dot x=g(x,\eps),
\end{equation}
with $g(x,\eps)$ given by
$$
\left(
\begin{array}{l}
x_2,\\
-x_1-\eps c x_2
+\dfrac{1-\alpha}{1-\eps\gamma(x_2-V)}+\alpha+\eps\beta(x_2-V)^2
\end{array}\right),
$$
possesses a closed orbit $\tilde x$ with the initial condition on $L$ (dotted curve in Fig~\ref{nfig}). The cycle we obtained is not a stick-slip limit cycle and therefore we want to rule out the possibility of the existence of $\tilde x.$ This will be accomplished over the Criterion of Bendixson \cite[Ch. X, \S7]{lef} which requires that the divergence of the vector field is strictly positive inside the limit cycle and along the limit cycle.
The divergence of (\ref{ns}) computes as 
$$
\begin{array}{l}   {\rm div} g(x)=\dfrac{\partial g_1}{\partial x_1}(x)+\dfrac{\partial g_2}{\partial x_2}(x)=\\
\qquad=-\eps c+\dfrac{\eps\gamma(1-\alpha)}{(1-\eps\gamma(x_2-V))^2}+2\eps\beta (x_2-V).\end{array}
$$
A feasible sufficient condition for ${\rm div} g(x)$ to be positive along $\tilde x$ is 
$$
  -c+\gamma(1-\alpha)+2\beta(\tilde x_2(t)-V)>0 \quad{\rm for\ all\ }t\in\mathbb{R},
$$
which needs more than just (\ref{abc}).

\vskip0.2cm

\noindent {\boldmath $x_2(\tau)>1:$} In his case the existence of a (finite-time stable) stick-slip limit cycle follows by applying Proposition~\ref{ppp}.

\section{Conclusion}  In this paper we used perturbation theory to obtain a sufficient condition (\ref{condi}) for the existence of stick-slip oscillations in a dry-friction oscillator on a moving belt (see (\ref{dfo}) and Fig.~\ref{mass_belt}) assuming that both the viscous friction and the nonlinear part of the dry friction characteristics are $\eps$-small. It can be seen from the proof of proposition~\ref{pr2} that our sufficient condition is sharp in the sense that the reversed inequality  (\ref{condi}) implies that no stick-slip limit cycles occur in (\ref{dfo}) as $\eps>0$ crosses zero. 
The test can be applied to virtually any nonlinear friction characteristics as long as it is small.
To illustrate the result, we used our test to derive a simple algebraic relation (\ref{abc}) for the constants of Stribeck friction characteristics that ensures the occurrence of stick-slip oscillations.

\section*{Acknowledgments} The research is supported by NSF Grant
CMMI-1436856.

\bibliographystyle{plain}


\end{document}